\documentclass{article}
\usepackage{amssymb}

\newtheorem{theorem}{Theorem}
\newtheorem{dfn}[theorem]{Definition}
\newtheorem{lemma}[theorem]{Lemma}
\newtheorem{proposition}[theorem]{Proposition}
\newtheorem{corollary}[theorem]{Corollary}
\newtheorem{remark}[theorem]{Remark}
\newtheorem{example}{Example}
\newenvironment{Proof}[1]{\paragraph{\sl Proof#1}}{$\square$}
\newcommand{\CC}{{\bf C}}
\newcommand{\RR}{{\bf R}}

\newcommand{\NN}{{\bf N}}

\newcommand{\kk}{{\bf k}}

\newcommand{\un}{\underline}

\newcommand{\cD}{{\mathcal D}}

\newcommand{\cH}{{\mathcal H}}

\newcommand{\cK}{{\mathcal K}}

\newcommand{\cN}{{\mathcal N}}

\newcommand{\Exp}{{\rm Exp}}

\def\mathclass#1#2{\noindent 2000 {\it Mathematics Subject
Classification}\/:\,\, {\sc Primary}\,#1;\,\, {\sc Secondary}\,#2}

\title{Gr\"obner $\delta$-bases and Gr\"obner bases \\
for differential operators\thanks{Research  supported by  DGESIC
BFM2001-3164, A.I.H.F.2000-0044 and FQM-218.}}
\date{October 30, 2001}
\author{Castro-Jim\'{e}nez,  F.J.\footnote{Facultad de Matem\'{a}ticas,
Universidad de Sevilla,  Sevilla, Spain. {E-mail: castro@cica.es}} \\
Moreno-Fr\'{\i}as, M.A.\footnote{Facultad de Ciencias, Universidad de
C\'{a}diz, C\'{a}diz, Spain. E-mail: mariangeles.moreno@uca.es}}

\begin{document}
\maketitle

\mathclass{13N10}{13P10,16S32,32C38,68W30.}

\abstract{This paper deals with the notion of Gr\"obner $\delta$-base
for some rings of linear differential operators by adapting the works
of W. Trinks, A. Assi, M. Insa and F. Pauer. We compare this notion
with the one of Gr\"obner base for such rings. As an application, and
following a previous work of A. Assi, we give some results on
finiteness and on flatness of finitely generated left modules over
these rings.}

\section{Introduction. }\label{introduction}
We will study Gr\"obner $\delta$-bases for some rings of linear
differential operators.

We have adapted to the differential case some notions and some
results obtained by  W. Trinks in \cite{trinks} and A. Assi (in
\cite{Assi} and \cite{assi1}) for the case of a commutative
polynomial ring with coefficients in a commutative unitary ring.

The notion of Gr\"obner $\delta$-base we introduce here is equivalent
to the one of Gr\"obner base defined by M. Insa and F. Pauer in
\cite{insa-pauer}. Nevertheless, we reserve the name Gr\"obner base
for the classical notion introduced in \cite{castro-thesis} (see also
\cite{castro2}). Besides the $\bf k$-algebras appearing in
\cite{insa-pauer}, the cases ${\mathcal H}={\bf k}[[X]][X^{-1}]$ and
${\mathcal H}={\bf k}\{X\}[X^{-1}]$ (when ${\bf k} = {\bf R}, {\bf
C}$) will be especially interesting in order to extend the results of
\cite{ACG1} and \cite{ACG2} to the rings of linear differential
operators with coefficients in ${\mathcal H}$.

Section 2 is devoted to the definition of the class of rings of
linear differential operators we will study and to the  theory of
Gr\"obner $\delta$-bases. We have, in these rings, a reduction
algorithm which allows the effective construction of a Gr\"obner
$\delta$-base for a given ideal, defined by a finite system of
generators. This is the aim of the sections 3,4 and 5.

In section 6 we compare the notions of Gr\"obner $\delta$-base and
Gr\"obner base in the case of the Weyl algebras. We prove that any
Gr\"obner base (in the sense of \cite{castro-thesis} (see also
\cite{castro2})) of a left ideal of a Weyl algebra is a Gr\"obner
$\delta$-base with respect to an appropriate well-ordering. We also
prove that the converse is not true.

We can deduce adapted algorithms for membership problem, elimination
problem and syzygies problem by using Gr\"obner $\delta$-bases
(instead of Gr\"obner bases) that could be in some cases with better
complexity.

In section 7 we apply previous results to the study of flatness of
some modules in a (local) relative situation. We also give a
finiteness results for some modules. These flatness results could be
compared to those of \cite{sabbah-1} for Rees modules over Rees
rings.

It is a pleasure to thank Professor A. Assi for his help and useful
suggestions.

\section{Gr\"obner $\delta$-bases. } \label{halgebra}

Here, $\kk$ is a field of zero characteristic. Let us denote by
$\kk[[X]]=\kk[[x_1,\ldots,x_n]]$ the ring of formal power series and
by $\kk((X))$ its quotient field.

Let us denote by $\kk((X))[\partial]=\kk((X))[\partial_1,\cdots,
\partial_n]$ the ring of linear differential operators with
coefficients in $\kk((X))$, where $\partial_i$ stands for the partial
derivative with respect to the variable $x_i$.

Let us consider a noetherian sub-${\kk}$-algebra ${\mathcal H}
\subset \kk((X))$ stable under the action of the partial derivatives
$\partial_1,\cdots,
\partial_n$. Let us denote by ${\mathcal D}$ (or ${\cal H}[\partial]$) the
sub-$\kk$-algebra (of $\kk((X))[\partial]$) of linear differential
operators generated by ${\mathcal  H}$ and $\{\partial_1, \cdots,
\partial_n \}$.

More generally, we will consider differential rings as
$\cD=\cH[\partial_1,\ldots,\partial_n]$ for any noetherian
sub-${\kk}$-algebra $\cH$ of
$\kk((\tilde{X}))=\kk((x_1,\ldots,x_n,x_{n+1},\ldots,x_{n+m}))$,
stable under the action of $\partial_i$ for $i=1,\ldots,n$.

The ring  ${\mathcal D}$ is the set of formal finite sums
$$\sum_{\alpha \in \NN^n}p_{\alpha}
\partial^{\alpha},
$$ where  $p_{\alpha}\in {\cH}$.

Let $<$ be a well-ordering, compatible with the sum in $\NN^n$ (i.e.
a well-ordering such that, for all $\gamma \in \NN^n$, we have
$\alpha + \gamma < \beta+\gamma $ if and only if $\alpha < \beta$).

\begin{dfn} \label{newton-diagram}
 Let  $P=\sum_{\alpha \in \NN^n } p_{\alpha} \partial^{\alpha}$
  be a non-zero  element of ${\mathcal D}$. The Newton $\delta$-diagram  of $P$ is  the set $$
  {\cN}^{\delta}(P)=
  \left \{
  \alpha \in \NN^{n}: p_{\alpha}\neq 0
  \right \}.
$$
\end{dfn}

\begin{dfn} \label{delta-exponent}
Let $P$ be a non-zero  element of  ${\mathcal D}$. We call the
element of\,\, $\NN^n$, $\mbox{max}_<\{{\cN}^{\delta}(P)\}$, the
$\delta$-exponent of $P$ with respect to $<$. It will be denote by
$\exp_<^{\delta}(P)$ or by $\exp^{\delta}(P)$ when no  confusion is
possible.
\end{dfn}

\begin{dfn} \label{delta-coefficient}
Let $P$ be a non-zero  element of  ${\mathcal D}$. We call the
element
 $p_{\alpha}\in {\cH}$, where $\alpha=\exp^\delta(P)$, the  $\delta$-coefficient of
 $P$ with respect to $<$. It will be denote by  $c_{<}^{\delta}(P)$
 or by $c^{\delta}(P)$ when no confusion is possible.
\end{dfn}

With these notations we have the following (see \cite[page
106-108]{moreno-tesis}):
\begin{lemma}\label{proex}
 Given two non-zero elements $P,Q$ in  ${\mathcal D}$. Then the following properties hold:
 \begin{enumerate}
  \item $\exp^\delta(PQ)=\exp^\delta(P)+\exp^\delta(Q)$ and
        $\exp^{\delta}([P,Q])< \exp^{\delta}(PQ)$.
  \item If $\exp^{\delta}(P)\neq \exp^{\delta}(Q) $ then
  $\exp^\delta(P+Q) = \mbox{max}_{<}\left \{
         \exp^\delta(P), \exp^\delta(Q) \right \}$.
  \item If $\exp^{\delta}(P)= \exp^{\delta}(Q) $ and $c^{\delta}(P)+c^{\delta}(Q)\neq 0$
        then  $\exp^\delta(P+Q) = \exp^\delta(P)=
        \exp^\delta(Q)$ and $c^{\delta}\left ( P+Q \right)=
        c^{\delta}(P)+c^{\delta}(Q)$.
 \item If $\exp^{\delta}(P)= \exp^{\delta}(Q) $ and $c^{\delta}(P)+c^{\delta}(Q)= 0$
        then  $\exp^\delta(P+Q) < \exp^\delta(P)$.
 \end{enumerate}
 \end{lemma}

All the ideals we will consider in ${\mathcal D}$ will be left
ideals.

Let $I$ be a non-zero ideal of ${\mathcal D}$. We denote by $$
\Exp^{\delta}_{<}(I)=\left\{ \exp^{\delta}(P): P \in I\setminus
  \{0\}\right\}\subseteq \NN^n.
  $$
  We write  $\Exp^{\delta}(I)$ when no
confusion is possible.

\begin{remark}\label{stair}
By Lemma \ref{proex} we have
$\Exp^{\delta}(I)+\NN^n=\Exp^{\delta}(I)$. So, by Dickson's Lemma
(see for example \cite{cox}), there is a finite generating subset $F$
of $\Exp^{\delta}(I)$, i.e. $$ \Exp^{\delta}(I)=\cup_{\alpha\in
F}(\alpha+\NN^n).$$ Any of the subset $F$ is called a $\delta$-stair
of $I$. \end{remark}

We denote by ${\cH}[\zeta]={\cH}[\zeta_1,\cdots, \zeta_n]$ the
(commutative) polynomial ring with coefficients in ${\cH}$ and with
variables $ \zeta_1, \cdots, \zeta_n $.

\begin{dfn} \label{delta-initial-form}
Let $P=\sum_{\alpha} p_{\alpha}\partial^{\alpha}$ be a  non-zero
element  of ${\mathcal D}$. The  $\delta$-initial form of $P$, with
respect to $<$, is
  $
     in_{<}^\delta(P)=c^{\delta}(P)\zeta^{\exp^{\delta}(P)} \in \cH[\zeta].
  $
   We write $in^\delta(P)$ when no confusion is possible.
  \end{dfn}

\begin{dfn} \label{delta-initial-ideal}
 Let $I$ be a non-zero ideal of  ${\mathcal D}$. The ideal generated by
 $$
 \left \{ in^{\delta}(P): P\in I\setminus \{0\}\right\}
 $$
  is called the
  $\delta$-initial ideal of $I$ with respect to $<$ (and it is denoted  by
 $in_{<}^\delta(I)$). We write $in^\delta(I)$ when no confusion is possible.
 \end{dfn}

\begin{remark} Note that $in^\delta(I)$ is a
$\zeta$-monomial ideal in ${\cH}[\zeta]$. If $I$ is generated by
$\{P_{1}, \cdots, P_{m}\}$ (i.e. $I=\cD(P_{1}, \cdots, P_{m})$), then
the ideals $ \cH[\zeta](in^{\delta}(P_1), \cdots,in^{\delta}(P_m))$
and $in^{\delta}(I)$ may be different.
\end{remark}

\begin{dfn} \label{grobner-delta-base}
Let $I$ be a non-zero ideal of  ${\mathcal D}$. A finite family
$\left \{ P_1, \cdots, P_m \right \}\subset I$ is called a Gr\"obner
$\delta$-base of $I$, with respect to the well ordering $<$, if $$
in^{\delta}(I)={\cH}[\zeta](in^{\delta}(P_1), \cdots,
in^{\delta}(P_m)).$$
\end{dfn}

\begin{remark} When $\cH=\kk$ (i.e.
in the ring $\kk[\partial]=\kk[\partial_1,\cdots, \partial_n]$), the
notion of Gr\"obner $\delta$-base and the one of Gr\"obner base
coincide. Here $\kk[\partial]$ is the ring of linear differential
operators with constant coefficients which is a commutative
polynomial ring.
\end{remark}

\section{Reduction in  ${\mathcal D}$. }\label{reduction}

Let $F$ be a  non-empty set in ${\mathcal D}$ and $\alpha \in
\NN^{n}$. Here we will use some notations of \cite{Assi}. Let
 $$
   K(\alpha;F)=\left \{c^{\delta}(P)\, : \,
   P\in F,\,  \alpha \in \exp^{\delta}(P)+\NN^{n}
  \right \}.
 $$

 We denote by  $C(\alpha;F)$ the ideal, in ${\cH}$, generated by
 $K(\alpha;F)$, i.e.
 $$
   C(\alpha;F)= {\cH}K(\alpha;F).
 $$
 If $K(\alpha;F)=\emptyset$ then  $C(\alpha; F)=\{0\}$.

\begin{example}
Let consider $F=\left\{ P_1,P_2 \right\}\subseteq
\kk[x_1,x_2][\partial_1,\partial_2]$, where
$P_1=a(x_1)\partial^2_1+b(x_2)\partial_1$ and
$P_2=c(x_2)\partial^2_2+e(x_1)\partial_2$ with $a(x_1),e(x_1)\in
\kk[x_1]$,  $a(x_1)\neq 0$ and $b(x_2),c(x_2)\in \kk[x_2]$ with
$c(x_2)\neq 0$. We consider the lexicographic order(c.f. \cite{cox})
in $\NN^{2}$ with $\partial_1
>\partial_2$.

Then  $\exp^\delta(P_1)=(2,0)$ and $\exp^\delta(P_2)=(0,2)$. Since
$K((1,1);F)=\emptyset$ we have
  $C((1,1);F)=\{0\}$. We also have
$$
  C\left ( (2,2);F\right)=\langle a(x_1),c(x_2)\rangle,
  \quad
  C\left ( (2,0);F\right)=\langle a(x_1)\rangle,
  \quad
  C\left ( (0,2);F\right)=\langle c(x_2)\rangle,
$$ where $\langle N \rangle$ stands for the ideal (in $\kk[x_1,x_2]$)
generated by $N$.
\end{example}

It is easy to check that
 \begin{enumerate}
  \item If $\beta \in \alpha + \NN^n$ then $K(\alpha; F)\subseteq K(\beta;F)$
  and $C(\alpha; F)\subseteq C(\beta; F)$.
  \item  If  $F_1 \subseteq F_2$ then  $K(\alpha; F_1)\subseteq K(\alpha; F_2)$
and  $C(\alpha; F_1)\subseteq C(\alpha; F_2)$ for all  $\alpha \in
\NN^n$.
\end{enumerate}

\begin{remark}\label{c(0,i)}
  If $I$ is a non-zero ideal of ${\mathcal D}$ then: \\ a)
  $
  K(\alpha; I)=\left\{c^{\delta}(P)\, : \, P\in I, \, \exp^{\delta}(P)=\alpha
  \right\}.
  $
Moreover, $K(\alpha; I)\bigcup \{0\}=C(\alpha; I)$. \\ b)
$C(0;I)=I\cap  {\cH}$.
\end{remark}

\begin{remark}\label{additional}
>From now on, we suppose ${\cH}$ verifying  two additional
conditions:
\begin{enumerate}
\item[1)] For any  subset $\{f_1,\cdots, f_r\}\subset {\cH}$ and for any   $f\in
{\cH}$ we can decide if $f\in \cH(f_1,\cdots, f_r)$, and in this
case, it is  possible to find  $q_1,\cdots, q_r \in {\cH}$ such that
$f=\sum_{i=1}^rq_if_i$.
\item[2)] For any subset $\{f_1,\cdots,f_r\}\subset {\cH}$ it is possible to find
a system of generators of the ${\cH}$-module of syzygies of
$\{f_1,\cdots, f_r\}$.
\end{enumerate}
\end{remark}

The algebras  $$ {\cH}=\kk[X],\,\, \kk(X),\,\, \kk[[X]],\,\,
\kk((X)), \,\, \kk[[X]][x_1^{-1},\,\, \cdots, x_n^{-1}]$$  and the
algebra $\kk\{X\}[x_1^{-1}, \cdots, x_n^{-1}]$ with $\kk=\RR $ or
$\CC$ verify  conditions 1) and 2).

\begin{dfn}
Let $F=\{P_1,\cdots, P_m\}\subseteq {\mathcal D} \setminus \{0\}$
with $P_i\neq 0$, $1\leq i \leq m$ and $P\in {\mathcal D}$. We will
say that
 $P$ is reduced with respect to $F$, if
  $\exp^{\delta}(P)\notin \bigcup_{i=1}^m\left(\exp^{\delta}(P_i)+\NN^n
 \right)
 $ or if  $\exp^{\delta}(P)\in \bigcup_{i=1}^m\left(\exp^{\delta}(P_i)+\NN^n
 \right)$ then $c^{\delta}(P)\notin C(\exp^{\delta}(P);F)$.
\end{dfn}

Let $F$ be a  non-empty subset of ${\mathcal D}\setminus \{0\}$. We
denote
  $$
  R(F)=\left\{R \in {\mathcal D}:\, R \mbox{ is reduced to  respect }F
  \right\}.
  $$

\begin{remark}
$R(F)$ is not necessarily a vector space over $\kk$.
\end{remark}

\begin{theorem}{\bf (Reduction algorithm)}\label{thdr}
Let $F=\{P_1, \cdots, P_m\} \subseteq {\mathcal D}$, with $P_i\neq
0$, $ i=1,\cdots, m$  and $P \in {\mathcal D}$. Then there exist $
Q_1, \cdots, Q_m,R \in {\mathcal D}$ such that
\begin{enumerate}
  \item $P=\sum_{i=1}^{m}Q_iP_i + R$.
  \item $R\in R(F)$.
  \item $\mbox{max}_{1\leq i \leq m}\{\exp^{\delta}(Q_iP_i), \exp^{\delta}(R)\}= \exp^{\delta}(P).$
\end{enumerate}
\end{theorem}

\begin{Proof}
We proceed by induction on $\exp^{\delta}(P)=\alpha$.\\ If $\alpha=
0$, then $P \in {\cH}$. So, we can consider two cases:
\begin{enumerate}
 \item If $0 \notin \bigcup_{i=1}^m\left(\exp^{\delta}(P_i)+\NN^n \right)$, then
 $$ P=\sum_{i=1}^m0P_i+P, \qquad \mbox{with}\qquad
P\in R(F).
 $$
 \item If $0 \in \bigcup_{i=1}^m\left(\exp^{\delta}(P_i)+\NN^n \right)$,
  then we consider the set
 $$
 \Lambda=\left\{i:\,  0 \in \exp^{\delta}(P_i)+\NN^n
 \right\} = \{i: \exp^\delta(P_i)=0\}.
 $$
  Thus for  $i\in \Lambda$ we have  $P_i \in {\cH}$ and
 we can consider two cases:
 \begin{enumerate}
  \item If $P \in \cH(P_i: i \in \Lambda)$ then
 $
   P=\displaystyle\sum_{i\in \Lambda}q_iP_i$ with $q_i\in {\cH}$ (according our assumption on $\cH$ we can calculate such
   elements $q_i$).
  In this case we have $P=\sum Q_iP_i +R$ where $Q_i=q_i$, for   $i\in \Lambda$; $Q_i=0$ for
 $i \notin \Lambda$ and $R=0$.

 \item If $P \notin \cH(P_i: i \in \Lambda)$ then $P\in
 R(F)$.
 \end{enumerate}
 \end{enumerate}

Suppose $\alpha > 0$ and the theorem proved  for
$\exp^{\delta}(P)<\alpha$.

Let $P\in {\mathcal D}$ be such that   $\exp^{\delta}(P)=\alpha$. We
have two possible cases:

 \begin{enumerate}
  \item If $\alpha \notin \bigcup_{i=1}^m\left(\exp^{\delta}(P_i)+\NN^n
  \right)$, then  $P=\displaystyle\sum_{i=1}^m0P_i+P$  and
  $P\in R(F)$.
  \item If $\alpha \in \bigcup_{i=1}^m\left(\exp^{\delta}(P_i)+\NN^n
  \right)$ then we consider the set
  $$
  \Lambda=\left\{i:\,  \alpha \in
  \exp^{\delta}(P_i)+\NN^n \right\}
  $$
  and the following two cases are possible:

  \begin{enumerate}
   \item If $c^{\delta}(P)\in C(\alpha;F)$,  then there exists
  $(q_i)_{i\in \Lambda} \in {\cH}$ such that
  $$
    c^{\delta}(P)=\sum_{i\in \Lambda}q_ic^{\delta}(P_i).
  $$
  We may write,
  $$
    P^{(1)}=P-\sum_{i\in \Lambda}q_i \partial^{\gamma^i}P_i,
    \qquad \mbox{ with } \qquad  \gamma^i+\exp^{\delta}(P_i)=\alpha.
  $$

  By construction,   $\exp^{\delta}(P^{(1)})<\exp^{\delta}(P)$.
  Hence, by induction, we may write
  $
    P^{(1)}=\sum_{i=1}^{m}Q'_iP_i + R'$,  with
  $R'\in R(F)$
 and finally
  $
    P=\sum_{i\notin \Lambda} Q'_iP_i+ \sum_{i\in \Lambda} \left(Q'_i+q_i
    \partial^{\gamma^i} \right)P_i+R'.$

 \item If $c^{\delta}(P)\notin C(\alpha; F)$ then  $P \in R(F)$.
 \end{enumerate}
\end{enumerate}

So, we have proved the existence of $Q_1,\ldots,Q_m,R$ verifying
conditions 1. and 2. of the statement. The condition 3. is easy to
verify. That ends the proof.
\end{Proof}

\begin{remark}
 We  call  $R \in {\mathcal D}$ a remainder of the reduction of
 $P$   by $\left ( P_1, \cdots, P_m \right )\subseteq {{\mathcal D}}^m$.
 We denote by
 $
 \widetilde{R}\left (P; P_1, \cdots, P_m \right ),
 $
 the set of remainders of the reduction of
  $P$ by  $\left \{ P_1, \cdots, P_m
\right \}$.

\end{remark}

\begin{remark} The proof of   Theorem \ref{thdr} provides an
algorithm to reduce an operator $P\in {\mathcal D}$ to respect  a
subset $F$ of ${\mathcal D}$.
\end{remark}

\begin{theorem}\label{tcdbR}
  Let $I$ be a non-zero ideal of ${{\mathcal D}}$ and $\left \{P_1, \cdots, P_r \right\} \subset I$.
  Then the following statements are equivalent:
 \begin{enumerate}
  \item  $\left \{P_1, \cdots, P_r \right\}$ is a Gr\"obner $\delta$-base of $I$.
  \item  For  $\alpha \in \NN^{n}$, we have
         $
     C(\alpha; I)=C(\alpha;P_1,\cdots, P_r).
         $
  \item For   $P\in I$ we have  $\widetilde {R}\left(P; P_1, \cdots,
        P_r \right)=\{0\}.$
 \end{enumerate}
\end{theorem}

\begin{Proof}{.}
 $1.\Longrightarrow 2.$:
  $C(\alpha;P_1, \cdots, P_r)$ is clearly contained in $C(\alpha;I)$. Conversely,  let
   $p(x)\in C(\alpha;I)$ then  there exists $P\in I\setminus \{0\}$
   such that   $in^{\delta}(P)=p(x)\zeta^{\alpha}$ and
   $in^{\delta}(P) \in in^{\delta}(I)$.  But by hypothesis, we have
  $$
 in^{\delta}(I)= {\cH}[\zeta]( in^{\delta}(P_1), \cdots, in^{\delta}(P_r)).
  $$
  Let us denote
   $$
   in^{\delta}(P_i)=p_i(x)\zeta^{\alpha_i}\qquad \mbox{with} \qquad 1 \leq i \leq
   r,$$
then
  $$
   p(x)\zeta^{\alpha}=\sum_{i=1}^rq_i(x,\zeta)p_i(x)\zeta^{\alpha_i},
  $$
  where
  $$
    q_i(x,\zeta)=\sum_\beta q_{i_{\beta}}(x)\zeta^{\beta} \in {\cH}[\zeta]\qquad {with}
  \qquad
   q_{i_{\beta}}(x)\in {\cH}.
  $$
  Thus,
  $$
  p(x)\zeta^{\alpha}
  =\sum_{i,\beta} q_{i_{\beta}}(x)p_i(x)\zeta^{\beta+\alpha_i},
  $$
 hence,
$$
  p(x)\zeta^{\alpha}\in {\cH}[\zeta]( p_i(x)\zeta^{\alpha_i}: \alpha\in \alpha_i+
  \NN^{n})
$$ and so, $$
  p(x)\in
 C(\alpha; P_1, \cdots, P_r).
 $$
Therefore $C \left (\alpha; I \right) \subseteq C \left ( \alpha;
P_1, \cdots, P_r \right )$ and it follows that
 $$
C\left (\alpha; I \right) = C \left ( \alpha; P_1, \cdots, P_r
\right). $$ \noindent $2.\Longrightarrow 3.$: Let $P\in I\setminus
\{0\}$, then by Theorem \ref{thdr}, there exists $Q_1, \cdots, Q_r,R
\in {\mathcal D}$ such that $$
 P=\sum_{i=1}^r Q_i P_i +R,
$$ where   $ R  \in \widetilde{R} \left ( P; P_1, \cdots, P_r \right
)$.

 Suppose  $R\neq 0$. Since
 $
 R=P-\sum_{i=1}^rQ_iP_i \in I,
 $
  we can consider two cases:
 \begin{enumerate}
  \item[i)] If  $\exp^{\delta}(R)\notin
  \bigcup_{i=1}^r\left(\exp^{\delta}(P_i)+\NN^n\right)$, then
   $
   C(\exp^{\delta}(R);P_1, \cdots,P_r)=(0).
   $
   Therefore
  $c^{\delta}(R)\notin  C(\exp^{\delta}(R);P_1, \cdots,P_r)$ and
   by hypothesis 2, $c^{\delta}(R)\notin C(\exp^{\delta}(R); I)$. But this is impossible
   since $R\in I$.
    \item[ii)] If $\exp^{\delta}(R)\in \bigcup_{i=1}^r\left(\exp^{\delta}(P_i)+\NN^n
  \right)$, then $$c^{\delta}(R)\notin  C(\exp^{\delta}(R);P_1,\ldots,P_r)=C(\exp^{\delta}(R);I)$$ because $R$
  is reduced with respect to  $\{P_1,\ldots,P_r\}$, and  this contradicts
that $R\in I$.
  \end{enumerate}
  Therefore $R=0$.

\noindent $ 3. \Longrightarrow 1.$: We must show that
 $
  in^{\delta}(I)={\cH}[\zeta]  (in^{\delta}(P_1), \cdots, in^{\delta}(P_r)).
$
Clearly $
 \langle in^{\delta}(P_1), \cdots,
in^{\delta}(P_r)
 \rangle \subseteq in^{\delta}(I)
$. Let
 $P\in I$ be a non-zero operator. Then we can write

$$P=p_{\alpha_0} \partial^{\alpha_0}+\widehat{P} $$ where
$p_{\alpha_0}\in {\cH}\setminus \{0\}$ and
$\exp^{\delta}(\widehat{P})<\alpha_0$.

 Then, by hypothesis and by Theorem \ref{thdr}, we have:
 $$
\alpha_0\in \bigcup_{i=1}^r\left(\exp^{\delta}(P_i)+\NN^n \right)
\qquad \mbox{and} \qquad c^{\delta}(P)\in C\left(\alpha_0;
P_1,\cdots, P_r \right).
 $$

   We consider the set $
   \Lambda=\left\{i\, : \,  \alpha_0\in \exp^{\delta}(P_i)+\NN^n\right\}.
   $ Then
  $
  c^{\delta}(P)=\sum_{i\in \Lambda}q_i^{(1)}c^{\delta}(P_i).
  $

 Let
 $$
 P^{(1)}=P-\sum_{i\in \Lambda}q_i^{(1)} \partial^{\gamma^i}P_i,
 \qquad
 \mbox{with}
 \qquad
 \gamma^i+\exp^{\delta}(P_i)=\alpha_0,
 $$
 then $P^{(1)}\in I$ and
 $\exp^{\delta}(P^{(1)})<\alpha_0$.\\

 Now we can  consider two cases:
 \begin{enumerate}
  \item[i)] If $P^{(1)}=0$, then
  $P=\sum_{i\in \Lambda}q_i^{(1)} \partial^{\gamma^i}P_i$
 and it can be checked that
  $$
  in^{\delta}(P)=\sum_{i \in \Lambda}q_i^{(1)}\zeta^{\gamma^i}in^{\delta}(P_i).
  $$
  \item[ii)] If $P^{(1)}\neq 0$, then by repeating the same procedure, we can obtain a family $P^{(k)}\in I$ with
  $\exp^{\delta}(P^{(k)})<\exp^{\delta}(P^{(k-1)})$. So, as $<$ is a well-ordering in $\NN^n$, there exists $l$,
  such that $P^{(l)}=0$.
 \end{enumerate}
 This completes the proof.
\end{Proof}

As a straightforward consequence of Theorem \ref{tcdbR} we get the
following result:

\begin{corollary}\label{exp-de-una-delta-base} Any Gr\"obner $\delta$-base of an ideal  $I\subseteq {\mathcal
D}$ is a system of generators of $I$. Moreover, if
$\{P_1,\ldots,P_r\}$ is a Gr\"obner $\delta$--base of $I$ then
$$\Exp^\delta(I)=\bigcup_{i=1}^r (\exp^\delta(P_i)+\NN^n).$$
\end{corollary}

\section{$S^{\delta}$-operators. }\label{s-sdelta-ope}

Let $F=\{P_1, \cdots, P_r\}\subseteq {\mathcal D}\setminus \{0\}$.
Let
   $$
     K(F)=\left \{
        \alpha \in \NN^{n}: \exists N \subseteq F,\,
        \alpha=lcm\{\exp^{\delta}(P); P \in N\}
          \right \},
          $$
          where $lcm$ stands for {\it less common multiple}, and
\label{dsicigias}
 $$
   F_{\alpha}=\left \{
   \left ( \lambda_1, \cdots, \lambda_r
   \right )\in {\cH}^{r}: \sum_{k=1}^{r}\lambda_kc^{\delta}(P_k)=0
   \mbox{ where } \lambda_k=0 \mbox{ if }\alpha \notin \exp^{\delta}(P_k)+\NN^{n}
    \right\} \subseteq {\cH}^r.
 $$

   $F_{\alpha}$   is isomorphic to the
  ${\cH}$-module of syzygies of
   $$
    \left \{ c^{\delta}(P_k): \alpha \in \exp^{\delta}(P_k)+\NN^{n},\,
     1\leq  k  \leq r \right \}.
     $$

Since ${\cH}$ is a noetherian algebra then   $F_{\alpha}$ is finitely
generated (as a ${\cH}$-module). Let
  $\left\{\left ( \lambda^{\tau}_1, \cdots, \lambda^{\tau}_r
   \right )\right\}$, $1\leq \tau \leq r_{\alpha}$, be a system of generators
   of
 $F_{\alpha}$.

 \begin{dfn}\label{s-delta-oper} With the notations as above,
for $\tau=1, \cdots, r_{\alpha}$, the element
  $$
   S_{\alpha,\tau}^{\delta}=
   \sum_{k=1}^{r}\lambda_{k}^{\tau}\partial^{\alpha-\exp^{\delta}(P_k)}
   P_k
 $$
   will be  called  a  $S^{\delta}$-operator of the
  set $F_{\alpha}$.
  \end{dfn}

\begin{proposition}\label{exp_deltapold}
  With the notations as above, we have
 $$
   \exp^{\delta}\left ( S_{\alpha,\tau}^{\delta}
 \right)< \alpha.
 $$
\end{proposition}

\begin{Proof}{.} We can write
 $$
 S_{\alpha,\tau}^{\delta}=
   \sum_{k=1}^{r}\lambda_{k}^{\tau}\partial^{\alpha-\exp^{\delta}(P_k)}
   P_k=
$$
 $$
  =\sum_{k=1}^{r}\lambda_{k}^{\tau}\partial^{\alpha-\exp^{\delta}(P_k)}
   c^{\delta}(P_k)\partial^{\exp^{\delta}(P_k)}+
 \sum_{k=1}^{r}\sum _{\beta < \exp^{\delta}(P_k)}\lambda_{k}^{\tau}\partial^{\alpha-\exp^{\delta}(P_k)}
 p_{\beta,k}\partial^{\beta}.
$$ Since
 $$
  \partial^{\alpha-\exp^{\delta}(P_k)}
   c^{\delta}(P_k)=c^{\delta}(P_k)\partial^{\alpha-\exp^{\delta}(P_k)}
  +A_{k}
  \qquad
  \mbox{with}
  \qquad  \exp^{\delta}(A_k) < \alpha-\exp^{\delta}(P_k),
  $$
$$
  \partial^{\alpha-\exp^{\delta}(P_k)}
   p_{\beta,k}=p_{\beta,k}\partial^{\alpha-\exp^{\delta}(P_k)}
  +B_{k}
   \qquad
  \mbox{with}
  \qquad  \exp^{\delta}(B_k) < \alpha-\exp^{\delta}(P_k)
  $$
  and
  $\sum_{k=1}^{r}\lambda_{k}^{\tau}c^{\delta}(P_k)=0,
  $
   finally
$$
 S_{\alpha,\tau}^{\delta}
  =
  \sum_{k=1}^{r}\lambda_{k}^{\tau}
  A_{k}\partial^{\exp^{\delta}(P_k)}+
 \sum_{k=1}^{r}\sum _{\beta < \exp^{\delta}(P_k)}\left(\lambda_{k}^{\tau}
 p_{\beta,k}\partial^{\alpha-\exp^{\delta}(P_k)+\beta}+
B_{k}\partial^{\beta}\right).
 $$
\end{Proof}

\begin{proposition}\label{tcD} Let $I$ be a non-zero ideal of  ${\mathcal
D}$ and
  $\left \{ P_1, \cdots, P_r \right \}$ be a system of generators of  $I$.
  Then the following  are equivalent:
 \begin{enumerate}
  \item $\left \{ P_1, \cdots, P_r \right \}$ is a  $\delta$-Gr\"obner base of
     $I$.
  \item For all  $P\in I$, we have
        $\widetilde{R}(P; P_1, \cdots, P_r)=\{0\}$.
  \item For all  $S^{\delta}$-operator, $S^{\delta}_{\alpha,\tau}$,
  of  $\left \{ P_1, \cdots, P_r \right \}$ we have
        $0\in \widetilde{R}\left (S^{\delta}_{\alpha,\tau};
         P_1, \cdots, P_r \right)$.
 \end{enumerate}
\end{proposition}

\begin{Proof}{.}

 $1. \Longrightarrow 2.$: See Theorem \ref{tcdbR}.\\

 $2. \Longrightarrow 3.$: Since  $S^{\delta}_{\alpha,\tau} \in I$,
 then, by assumption,
 $0\in \widetilde{R}\left (S^{\delta}_{\alpha,\tau};
         P_1, \cdots, P_r \right)$.\\

 $3. \Longrightarrow 1.$:
 Let $P \in I$ be a non-zero operator. We must show that  $in^{\delta}(P)\in
{\cH}[\zeta] (in^{\delta}(P_{1}), \cdots, in^{\delta}(P_{r})) $. We
may write
   $P=\sum_{i=1}^rH_iP_i$, with  $H_i \in {\cH}[\partial]$.

 Suppose $$\alpha_0=\mbox{ max}_{i}
 \left \{ \exp^{\delta}(H_iP_i)\right \}\qquad  \mbox{and}
 \qquad
   \exp^{\delta}\left ( H_{i_k}P_{i_k}\right)=\alpha_0,
   \qquad
   k=0, \cdots, t.
 $$
 Hence,  by Lemma \ref{proex},
 $$
   \exp^{\delta}\left ( H_{i_k}\right)+
  \exp^{\delta}\left ( P_{i_k}\right)=\alpha_0,
  \qquad
   k=0, \cdots, t.
 $$
 We can consider two cases:
 \begin{enumerate}
  \item[a)] If
   $
     \sum_{k=0}^{t}c^{\delta}\left ( H_{i_k}\right)
      c^{\delta}\left ( P_{i_k}\right)\neq 0
   $
   then
   $$
     in^{\delta}(P)=c^{\delta}(P)\zeta^{\alpha_0},
   \mbox{ with }
   c^{\delta}(P)=\sum_{k=0}^{t}c^{\delta}\left ( H_{i_k}\right)
      c^{\delta}\left ( P_{i_k}\right).
   $$
 Therefore,
 $$
   in^{\delta}(P)=\sum_{k=0}^{t}
c^{\delta}\left ( H_{i_k}\right)
      c^{\delta}\left ( P_{i_k}\right)\zeta^{\alpha_0}
 =\sum_{k=0}^{t}
c^{\delta}\left (
H_{i_k}\right)\zeta^{\alpha_0-\exp^{\delta}(P_{i_k})}
 in^{\delta}(P_{i_k})
 $$
 and so,
 $in^{\delta}(P)\in{\cH}[\zeta] (in^{\delta}(P_{1}), \cdots, in^{\delta}(P_{r}))$.

\item[b)] Suppose now $\sum_{k=0}^{t}c^{\delta}\left ( H_{i_{k}}\right)
c^{\delta}\left ( P_{i_{k}}\right)=0$. Let us denote
$\alpha^i=\exp^{\delta}(P_i)$, $i=1,\ldots,r$; we consider the set
$
 \Lambda=\left\{i:\,
 \exp^{\delta}(H_i)+\alpha^i=\alpha_0
 \right\},
 $ and we suppose  $\gamma=lcm \{ \exp^{\delta}(P_i):\, i\in \Lambda\}$.

 We may  write
  $$
   P=\sum_{i \notin \Lambda}H_iP_i+\sum_{i\in
   \Lambda}c^{\delta}(H_i)\partial^{\exp^{\delta}(H_i)}P_i+
   \sum_{i\in
   \Lambda}(H_i-c^{\delta}(H_i)\partial^{\exp^{\delta}(H_i)})P_i.
 $$
 We can  identify  $\left(c^{\delta}(H_i) \right)_{i \in
 \Lambda}$ with an element of $F_{\gamma}$. Let
 ${\un \lambda}^1, \cdots,\un{\lambda}^p $ be a family  of generators of
 $F_{\gamma}$ where
 $$
 \un{\lambda}^{\tau}=(\lambda_1^{\tau}, \cdots,
 \lambda_r^{\tau})\qquad
 \mbox{with}
 \qquad \lambda_j^{\tau}=0 \quad \mbox{if} \quad  \gamma \notin \exp^{\delta}(P_j)+\NN^n.$$
 Now for each $i \in \{1,\ldots, r\}$ we define $s_i$ as follows:
 $$
  s_i=\left\{ \begin{array}{lcr}
      c^{\delta}(H_i) & \mbox{if }& i\in \Lambda\\
               0 & \mbox{if }& i\notin \Lambda.
              \end{array}
  \right.
 $$
 Hence $\un s=(s_1, \cdots, s_r)\in F_\gamma$, and then there exist $u_1,\cdots, u_p \in
 {\cH}$ such that $\un{s} = \sum_{\tau=1}^p
 u_{\tau}\un{\lambda}^{\tau}$. Let us denote  $\beta^i=\exp^{\delta}(H_i)$,
 for  $i \in \Lambda$, then
 $$
 \sum_{i\in \Lambda}c^{\delta}(H_i)\partial^{\beta^i}P_i=
 \sum_{\tau=1}^p u_{\tau} \left( \sum_{i=1}^r
 \lambda_i^{\tau} \partial^{\beta^i}P_i \right).
 $$

The element $\alpha_0$ is, by definition, a common multiple of the
elements $\{\exp^{\delta}(P_i):\,  i \in \Lambda \}$ then  there
exists $\epsilon \in \NN^n$ such that  $\alpha_0=\gamma + \epsilon$
and so
 $\beta^i=\gamma-\alpha^i+\epsilon$.\\
  If $j \notin \Lambda$ and
$\gamma \in \exp^{\delta}(P_j)+\NN^n$ we denote  $\beta^j=
\gamma-\alpha^j + \epsilon$. Therefore,
 $$
\sum_{i\in \Lambda}c^{\delta}(H_i)\partial^{\beta^{i}}P_i
=\sum_{\tau=1}^p u_{\tau}\left( \sum_{i=1}^r \partial^{\epsilon}
\lambda_i^{\tau}\partial^{\gamma-\alpha^i}P_i\right)+ \sum_{\tau=1}^p
u_{\tau}\left( \sum_{i\vert \gamma-\alpha^i>0} B_i^{\tau}
\partial^{\gamma-\alpha^i}P_i\right)
$$ where $\exp^{\delta}(B_i^{\tau})<\epsilon$.\\ Therefore, by
Definition \ref{s-delta-oper},
 $$
 \sum_{i
\in \Lambda}c^{\delta}(H_i)\partial^{\beta^i}P_i= \sum_{\tau=1}^p
u_{\tau}\partial^{\epsilon}S_{\gamma,\tau}^{\delta}+ \sum_{i\vert
\gamma-\alpha^i>0}\left( \sum_{\tau=1}^p
u_{\tau}B_{i}^{\tau}\right)\partial^{\gamma-\alpha^i}P_i. $$ But by
hypothesis, we have
 $$
  S^{\delta}_{\gamma,\tau}=\sum_{j=1}^rQ_{j}^{\gamma,\tau}P_j,
$$ with $\gamma > \exp^{\delta}(S_{\gamma,\tau}^{\delta})=\mbox{ max
}_{1\leq j \leq r} \{\exp^{\delta}(Q_j^{\gamma,\tau}P_j)\}$. Hence,
 $$ \sum_{i\in
\Lambda}c^{\delta}(H_i)\partial^{\beta^i}P_i= \sum_{j=1}^r\left(
\sum_{\tau=1}^p
u_{\tau}\partial^{\epsilon}Q_{j}^{\gamma,\tau}\right)P_j+ \sum_{j
\vert \gamma-\alpha^j>0}\left(\sum_{\tau=1}^{p}u_{\tau}B_j^{\tau}
\right)\partial^{\gamma-\alpha^j}P_j. $$ Therefore
 $$
P=\sum_{i=1}^rH'_iP_i $$
 where
\begin{itemize}
\item If $i \in \Lambda$,
      $$
H'_i=H_i-c^{\delta}(H_i)\partial^{\beta^i}+
      \sum_{\tau=1}^pu_{\tau}\partial^{\epsilon}Q_i^{\gamma,\tau}+
      \sum_{\tau=1}^p u_{\tau}B_i^{\tau}\partial^{\gamma-\alpha^i}.
     $$
\item If $i \notin \Lambda$ and $\gamma-\alpha^i >0$,
     $$
       H'_i=H_i+\sum_{\tau=1}^p u_{\tau}\partial^{\epsilon}Q_i^{\gamma, \tau}+
       \sum_{\tau=1}^p u_{\tau}B_i^{\tau}\partial^{\gamma-\alpha^i}.
     $$
\item If $i \notin \Lambda$ and $\gamma-\alpha^i$ is not greater than 0,
     $$
     H'_i=H_i+\sum_{\tau=1}^p u_{\tau}\partial^{\epsilon}Q_i^{\gamma,\tau}.
     $$
\end{itemize}

Hence, we have obtained an expression for $P$ as a combination of the
$P_i$ where $\exp^{\delta}(H'_iP_i)< \alpha_0$, then
$\mbox{max}_i\{\exp^{\delta}(H'_iP_i)\}< \alpha_0$. But this process
stops because $<$ is a well-ordering in $\NN^n$. So, there exists an
expression of $P$ with the conditions of the case $a)$.
\end{enumerate}

\end{Proof}

\section{Construction of a Gr\"obner $\delta$-base. }\label{algoritmo}

Let $I$ be a non-zero ideal of ${\mathcal D}$ and let $F=\left\{
P_1,\ldots, P_r \right\}$ be a system of generators of $I$. We will
show here how to build a Gr\"obner  $\delta$-base of the ideal $I$
(with respect to a ordering $<$). We will follow the main lines of
Buchberger's algorithm, adapted to our case (see \cite{buchberger1},
\cite{trinks} and \cite{Assi}).

 Let
 $K(F)=\left\{\alpha^1,\ldots, \alpha^s
 \right\}$ (see Section \ref{s-sdelta-ope}). Let
 $\left\{S^{\delta}_{\alpha^j,\tau}\right\}$,$1 \leq j \leq s$, $1\leq \tau \leq
 r_j$ the family of $S^{\delta}$-operators associated to $F$.

 We suppose that $\left\{ P_1,\ldots, P_r \right\}$ is not a Gr\"obner
 $\delta$-base for $I$, then (by Proposition \ref{tcD}) there
 exists
$S^{\delta}_{\alpha_0,\tau}$ such that  $0\notin \widetilde{R}\left
(S^{\delta}_{\alpha_0,\tau}; P_1, \cdots, P_r \right),
$
then  let  $$ P_{r+1}\in \widetilde{R}\left
(S^{\delta}_{\alpha_0,\tau}; P_1, \cdots, P_r  \right) $$ and we
repeat this process with $\{P_1, \cdots, P_r, P_{r+1}\}$.
\begin{remark}
 If a $S^{\delta}$-operator, $S$, of $F$ verify that $0 \in \widetilde{R}\left
(S; P_1, \cdots, P_r \right)$ then $0 \in \widetilde{R}\left (S; P_1,
\cdots, P_r, P_{r+1} \right)$.
\end{remark}

The following Proposition assures that this procedure terminates

\begin{proposition}\label{algoritmodebase}  With the notations as
above, there exists
 $\rho \in \NN$ such that for all  $S^{\delta}$-operator $S$ of
 $\left \{
 P_1, \cdots, P_{r+\rho}
 \right\}$  we have  $0\in \widetilde{R}\left
(S;P_1, \cdots, P_{r+\rho} \right ). $
\end{proposition}
\begin{Proof}{.} See
\cite[pages 131-133]{moreno-tesis}.
\end{Proof}

\section{Gr\"obner bases and the
Gr\"obner $\delta$-bases. }\label{grobnerydeltabases} In this section
we will work on the Weyl algebra $A_n(\kk)=\kk[X][\partial]$, so we
suppose here $\cH=\kk[X]=\kk[x_1,\ldots,x_n]$.

Let $<_{x}$, $<_{\partial}$ be monomial orderings  in $\NN^n$.

We denote by $X^{\alpha}\partial^{\beta}$ the monomial
$$x_1^{\alpha_1}\cdots x_n^{\alpha_n}\partial_1^{\beta_1}\cdots
\partial_n^{\beta_n}$$  Let us define on $\NN^n \times \NN^n$ the total ordering
(denoted $<$) by

$$
  (\alpha(1), \beta(1))<
  (\alpha(2), \beta(2))
 \Longleftrightarrow
 \left \{ \begin{array}{ccc}
      \beta(1) & <_{\partial} & \beta (2)\\
      \mbox{or} & & \\
      \beta(1)=\beta(2) & \mbox{ and }  & \alpha(1) <_{x}
      \alpha(2).
          \end{array}
 \right.
$$

\begin{remark}
The relation $<$, defined in $\NN^n \times \NN^n$, is a monomial
ordering. This  well-ordering is called an elimination order (see for
example {\rm \cite{cox}}).
\end{remark}

For the notion of Gr\"obner base on $A_n(\kk)$ and some related
results we follow here \cite{castro-thesis} (see also
\cite{castro2}).

\begin{theorem}
Let $G=\left \{P_1, \cdots, P_r \right\}$ be a  system  of generators
for a non-zero ideal $I\subset{A_n(\kk)}$. Then if $G$ is a Gr\"obner
base for $I$, with respect to $<$, then $G$ is a Gr\"obner
$\delta$-base for $I$ with respect to $<_{\partial}$.
\end{theorem}

\begin{Proof}{.}
Let $P \in I$ be a non-zero operator. We  must  show  that
$$in^{\delta}(P)\in {\cH}[\zeta]( in^{\delta}(P_1), \cdots,
in^{\delta}(P_r)).$$ For $i=1,\ldots, r$, we may write $
P_i=a_i\partial^{\alpha_i}+ \widehat{P_i} $ where
$\exp^{\delta}(P_i)=\alpha_i$, $\exp^{\delta}(\widehat {P_i})<
\alpha_i$ and
 $a_i \in {\cH}$.
Thus
$
 in^{\delta}(P_i)=a_i\zeta^{\alpha_i}.
$ By the division algorithm in $A_n(\kk)$, (see \cite{castro-thesis}
and \cite{castro2}) there exists $Q_{i_1}, \cdots, Q_{i_N}\in
A_n(\kk)$, $1 \leq i_j \leq r$,
 satisfying
$  P=Q_{i_1}P_{i_1}+ \cdots +Q_{i_N}P_{i_N}$ where
$\exp_{<}(Q_iP_i)\neq \exp_{<}(Q_jP_j)$ for $i \neq j$.

We can suppose $$
  \exp_{<}(Q_{i_N}P_{i_N})< \exp_{<}(Q_{i_{N-1}}P_{i_{N-1}}) < \cdots
  <
  \exp_{<}(Q_{i_1}P_{i_1}).
$$
 We can write
 $$
  Q_{i_j}=c_{i_j}\partial^{\beta_{i_j}}+ \widehat{Q}_{i_j}
$$ where  $\exp^{\delta}(Q_{i_j})= \beta_{i_j}$,
$\exp^{\delta}(\widehat{Q}_{i_j})<\beta_{i_j}$, $c_{i_j}\in {\cH}$.
Thus
 $
  \exp_{<}(Q_{i_j})=(\exp_{<_{x}}(c_{i_j}),\beta_{i_j}).
$

Therefore, $$
  P=\sum_{j=1}^{N}
  c_{i_j} a_{i_j}\partial^{\beta_{i_j}+ \alpha_{i_j}}+
\sum_{j=1}^{N}c_{i_j}A_{i_j}\partial^{\alpha_{i_j}} + \sum_{j=1}^{N}
c_{i_j}\partial^{\beta_{i_j}}\widehat{P}_{i_j}+
\sum_{j=1}^{N}\widehat{Q}_{i_j}a_{i_j}\partial^{\alpha_{i_j}}+
\sum_{j=1}^{N}\widehat{Q}_{i_j}\widehat{P}_{i_j} $$ where,
  $\exp^{\delta}\left( \sum_{j=1}^{N}
  c_{i_j}a_{i_j}\partial^{\beta_{i_j}+ \alpha_{i_j}}\right )
\leq \mbox{max}_{1\leq j \leq N} \{\beta_{i_j}+\alpha_{i_j}\} $,

 $\exp^{\delta}\left( \sum_{j=1}^{N}
c_{i_j}A_{i_j}\partial^{\alpha_{i_j}}
  \right)<\mbox{max}_{1\leq j \leq N}
\{\beta_{i_j}+\alpha_{i_j}\} $,

 $\exp^{\delta}\left(
\sum_{j=1}^{N}c_{i_j}\partial^{\beta_{i_j}}\widehat{P}_{i_j} \right)
<\mbox{max}_{1\leq j \leq N} \{\beta_{i_j}+\alpha_{i_j}\}$,

 $\exp^{\delta}\left(\sum_{j=1}^{N}\widehat{Q}_{i_j}a_{i_j}
\partial^{\alpha_{i_j}}\right)
<\mbox{max}_{1\leq j \leq N} \{\beta_{i_j}+\alpha_{i_j}\} $ and

$\exp^{\delta}\left(\sum_{j=1}^{N}\widehat{Q}_{i_j}\widehat{P}_{i_j}\right)
<\mbox{max}_{1\leq j \leq N} \{\beta_{i_j}+\alpha_{i_j}\}.$

Let  $j_0$ be such that
 $$
  \beta_{i_{j_0+1}}+\alpha_{i_{j_0+1}}< \beta_{i_{j_0}}+\alpha_{i_{j_0}}=
  \beta_{i_{j_0-1}}+\alpha_{i_{j_0-1}}=\cdots=
\beta_{i_1}+\alpha_{i_1}.
  $$
  Since
  $
  \sum_{j=1}^{j_0} c_{i_j}a_{i_j} \neq 0,
  $
  we have
  $
  in^{\delta}(P)=\left ( \sum_{j=1}^{j_0} c_{i_j}a_{i_j}\right)
\zeta^{\beta_{i_1}+\alpha_{i_1}}.$ Therefore, $$
  in^{\delta}(P)=
\sum_{j=1}^{j_0}c_{i_j}in^{\delta}(P_{i_j})\zeta^{\beta_{i_j}} $$ and
so
$
in^{\delta}(P)
\in
{\cH}[\zeta]( in^{\delta}(P_1), \cdots, in^{\delta}(P_r)) $.  This
completes the proof.
\end{Proof}

The converse result is not true as we show in the following example:

\begin{example}\label{ejemplo-bases-y-delta-bases}
Let $I\subset A_2(\CC)={\bf C}[x_1,x_2][\partial_1,\partial_2]$ be
the left ideal generated by the operators
 $$
   P_1= x_1\partial_1+a\partial_2+b, \quad
   P_2=(x_2-x_1)\partial_2-d
 $$ with $a,b,d \in {\bf C}[x_1,x_2]$.

We will prove\footnote{By using the degree lexicographical order with
$\partial_2 <_\partial \partial_1$ and $x_2 <_x x_1$} that $\{P_1,P_2
\}$ is a Gr\"obner $\delta$-base which is not a Gr\"obner base of
$I$, for a particular choice of the polynomials $a,b,d$.

We have  $$ \exp_{<}(P_1)=(1,0,1,0), \qquad \exp_{<}(P_2)=(1,0,0,1).
$$ Then $$S(P_1, P_2)=\partial_2P_1+\partial_1P_2=
x_2\partial_1\partial_2+\partial_2a\partial_2+\partial_2b-\partial_1d-\partial_2,
$$ and then $$ \exp_{<}\left( S(P_1, P_2) \right)=(0,1,1,1)\notin
\left<
 (1,0,1,0), (1,0,0,1)
\right>= \left <
 \exp_{<}(P_1), \exp_{<}(P_2)
\right>. $$

So, $G=\{P_1,P_2\}$ is not a Gr\"obner base of the ideal $I$, for any
$a,b,d\in \CC[x_1,x_2]$.

We will prove that, for some $a,b,d\in \CC[x_1,x_2]$, the set
$G=\{P_1,P_2\}$ is a Gr\"obner $\delta$-base of $I$.

We have $$
   \exp^{\delta}(P_1)=(1,0), \qquad c^{\delta}(P_1)=x_1
$$ and $$
   \exp^{\delta}(P_2)=(0,1), \qquad c^{\delta}(P_2)=x_2-x_1.
$$

We will compute the associated $S^{\delta}$-operators (see Definition
\ref{s-delta-oper}).

As $$
   \alpha=lcm((1,0),(0,1))=(1,1)
$$ we must first compute  a system of generators of $$ F_{(1,1)}(P_1,
P_2)= \left \{ (\lambda_1,\lambda_2)\in {\bf C}[x_1,x_2]:
\lambda_1c^{\delta}(P_1)+\lambda_2c^{\delta}(P_2)=0 \right\}.$$ In
fact we have  $$ Syz(c^{\delta}(P_1), c^{\delta}(P_2))=F_{(1,1)}(P_1,
P_2)= \left <
 (x_2-x_1, -x_1)
\right > $$ and then  $$
  S^{\delta}_{(1,1),(x_2-x_1,-x_1)}=(x_2-x_1)\partial^{(1,1)-(1,0)}P_1-x_1\partial^{(1,1)-(0,1)}P_2=
$$ $$ =(x_2-x_1)a\partial_2^2+(x_2-x_1)\partial_2(a)\partial_2+
 (x_2-x_1)b\partial_2+
$$ $$ +(x_2-x_1)\partial_2(b)+
x_1\partial_2+x_1d\partial_1+x_1\partial_1(d).
 $$
Now we reduce $ S^{\delta}_{(1,1),(x_2-x_1,-x_1)}$ by  $(P_1,P_2),$
say $$ S^{\delta}_{(1,1),(x_2-x_1,-x_1)}-d P_1-\partial_2aP_2-bP_2=
 $$
$$
 =(x_2-x_1)\partial_2(b)+
x_1\partial_2+x_1\partial_1(d)-a\partial_2+a\partial_2(d)+\partial_2(a)d=
$$
 $$
=(x_2-x_1)\partial_2(b)+(x_1-a)\partial_2+x_1\partial_1(d)+\partial_2(ad).\\
$$

Then $\{P_1,P_2\}$ is a Gr\"obner $\delta$-base of $I$ if $a=x_1,\,
b\in {\bf C}[x_1]$ and $d\in {\bf C}$.
\end{example}

\begin{remark}\label{remark-after-sun} The example before proves a little more. Let us
consider $\cH=\CC[x_1,\ldots,x_n]$ (for $n\geq 3$) and the ring of
differential operators $\cD=\cH[\partial_1,\partial_2]$ (which is a
sub-algebra of the Weyl algebra $A_n(\CC)$).

Let $I\subset \cD$ be the left ideal generated by the operators
 $$
   P_1= x_1\partial_1+a\partial_2+b, \quad
   P_2=(x_2-x_1)\partial_2-d
 $$ with $a,b,d \in {\bf C}[x_1,\ldots,x_n]$. An analogous
 computation to the one of example \ref{ejemplo-bases-y-delta-bases}
 proves that $\{P_1,P_2\}$ is a Gr\"oner $\delta$-base of $I$ if
 $a=x_1, b\in \CC[x_1,x_3,\ldots,x_n]$ and $d\in \CC[x_3,\ldots,x_n].$
\end{remark}

\section{Applications: Flatness and finiteness. }

As elementary applications of  Gr\"obner $\delta$-bases, we have the
effective solution for the ideal membership problem, variable
elimination problem and effective intersection of ideals. We also can
calculate a generating system of the ${\cH}[\partial]$-module of
syzygies of a finite subset $\left\{P_1,\cdots,P_r \right\}$ of
${\cH}[\partial]$, as well as a free resolution of a finitely
generated (left) ${\cH}[\partial]$-module. Calculating free
resolutions of a ${\cH}[\partial]$-module,  we have found examples
where the use of Gr\"obner $\delta$-bases is, in some sense, more
efficient that the one of Gr\"obner bases (see \cite[page
189-190]{moreno-tesis}).

In this section the ring $\cH$ is a noetherian sub-${\kk}$-algebra of
$$\kk((\tilde{X}))=\kk((x_1,\ldots,x_n,x_{n+1},\ldots,x_{n+m})),$$
stable under the action of $\partial_i$ for $i=1,\ldots,n$ and
satisfying the two additional conditions of Remark \ref{additional}.
We denote as before
$\cD=\cH[\partial]=\cH[\partial_1,\ldots,\partial_n]$.

The aim of this section is to characterize flatness and finiteness of
a ${\mathcal D}$-module by using the notion of Gr\"obner
$\delta$-bases, following the work of A. Assi \cite{assi1} in the
commutative case.

We can see the quotient $\cD/I$ as a family of $A_n(\kk)$-modules,
the space of parameters being $\CC^{m}$. In this section we will see
when this family is flat.

Let $S$ be multiplicatively closed subset of ${\cH}$. The ring
$S^{-1}{\cH}$ is a noetherian sub-$\kk$-algebra of
$\kk((\tilde{X}))$, stable under the action of the derivations
$\partial_1,\cdots,
\partial_n$ and satisfying the two additional conditions of Remark \ref{additional}.
So, we can consider the sub-$\kk$-algebra $S^{-1}\cD$ of
$\kk((\tilde{X}))[\partial]$, generated by $S^{-1}\cH$ and
$\partial_1,\ldots,\partial_n$.

One can define in $S^{-1}\cD$ the notions of section \ref{halgebra}.

Let $I\subset \cD$ a left ideal. We denote by $S^{-1}I$ the ideal of
$S^{-1}\cD$ generated by $I$  and  by $in^{\delta}(S^{-1}I)$ the
ideal (of $S^{-1}{\cH}[\zeta]$) generated by $\left\{in^{\delta}(P):
P\in (S^{-1}I)\setminus \{0\} \right\}$. Here $S^{-1}{\cH}[\zeta]$
denotes the polynomial ring in the variables
$\zeta=(\zeta_1,\ldots,\zeta_n)$ and coefficients in $S^{-1}\cH$. We
have:

\begin{proposition}\label{base-de-sI}
Suppose $\left\{P_1,\cdots, P_r \right\}$ is a Gr\"obner
$\delta$-base of $I$. If $S$ be a multiplicatively closed subset of
${\cH}$ then $in^{\delta}(S^{-1}I)$ is generated by $\left\{
in^{\delta}(P_1),\cdots, in^{\delta}(P_r)\right\} $ in
$S^{-1}{\cH}[\zeta]$. In particular, $in^{\delta}(S^{-1}I) =
S^{-1}(in^\delta(I))$ and $\{P_1/1,\ldots,P_r/1\}$ is a Gr\"obner
$\delta$--base of $S^{-1}I$.
\end{proposition}

Let ${{\mathfrak P}}$ be a prime ideal of ${\cH}$. Then
$S={\cH}\setminus {{\mathfrak P}}$ is a multiplicatively closed
subset of ${\cH}$. We denote ${\cH}_{{\mathfrak P}}=S^{-1}{\cH}$,
$\cD_{\mathfrak P}= S^{-1}\cD$ and $I_{{\mathfrak P}}=S^{-1}I$ and .

For each ideal $\cK$ in $\cH$, we denote $ V(\cK)=\left\{{{\mathfrak
P}}\in Spec({\cH}): \, \cK \subseteq {{\mathfrak P}} \right\},$ which
is a Zariski closed subset of $Spec(\cH)$. Here we endowed the set
$Spec(\cH)$  of prime ideals of $\cH$ with its Zariski topology.

Let consider $J=\prod_{i=1}^s C\left(\alpha(i);I\right)$ as an ideal
in $\cH$, where $\left\{\alpha(1),\cdots, \alpha(s)\right\}$ is a
$\delta$-stair of the ideal $I$ (see Remark \ref{stair}). Let us
denote $U= Spec({\cH})\setminus V(J)$.  We have:

\begin{theorem}\label{flateness} With the notations as above, let ${{\mathfrak P}} \in
U$. Then ${\cD}_{{\mathfrak P}}/I_{{\mathfrak P}}$ is a free (and
then a flat) ${\cH}_{{\mathfrak P}}$-module.
\end{theorem}
\begin{Proof}{.}
Let $M$ be the free ${\cH}_{{\mathfrak P}}$-module generated by
$\left\{\partial^{\alpha}:\, \alpha \in \NN^n \setminus
\Exp^\delta(I) \right\}$. Obviously we have
$\Exp^\delta(I)=\Exp^\delta(I_{\mathfrak P})$. Let us consider a
Gr\"obner $\delta$-base $\{P_1,\ldots,P_r\}$ of $I$. By Proposition
\ref{base-de-sI}, $\{P_1/1,\ldots,P_r/1\}$ is a Gr\"obner
$\delta$-base of $I_{\mathfrak P}$. Now, applying the reduction
algorithm with respect to $\{P_1/1,\ldots,P_r/1\}$ (see Theorem
\ref{thdr}), each $P\in \cD_{\mathfrak P}$ can be written as a sum
$$P=P'+P''$$ with $P'\in I_{\mathfrak P}$ and $P''\in M$. Here we
have used the equality $C(\alpha(i);I_{\mathfrak P}=\cH_{\mathfrak
P}$ for each $i=1,\ldots,s$.

So, we have proved that ${\cD}_{{\mathfrak P}}=I_{{\mathfrak P}}+M$
and it is obvious that $I_{{\mathfrak P}}\cap M=(0)$, so the
$\cH_{\mathfrak P}$--modules ${\cD}_{{\mathfrak P}}/I_{{\mathfrak
P}}$ and $M$ are isomorphic. Then $M$ is a free ${\cH}_{{\mathfrak
P}}$-module.
\end{Proof}

\begin{proposition} With the notations as above, we have
\begin{enumerate}
 \item If  $C(0;I)=I \cap {\cH}\neq (0)$, then
 $U=Spec({\cH})\setminus V\left( C(0;I) \right)$ is the maximal open set of flatness.
  \item If $C(\alpha(k);I)={\cH}$ for each $k \in \{1,\ldots,s\}$, then ${\cD}/I$ is a
  flat ${\cH}$-module.
\end{enumerate}
\end{proposition}
\begin{Proof}{.}
\begin{enumerate} \item We have  $C(0;I)={\cH} \cap I$ (see Remark \ref{c(0,i)}).
Suppose $U$ is not  maximal, then there exists ${\mathfrak P}\in
Spec({\cH})\setminus U$ such that ${\cD}_{\mathfrak P}/I_{\mathfrak
P}$ is ${\cH}_{\mathfrak P}$-flat. If $C(0;I)\neq (0)$ then
${\cH}_{\mathfrak P} \cap I_{\mathfrak P}\neq (0)$, which is
impossible by flatness of ${\cH}_{\mathfrak P}[\partial]/I_{\mathfrak
P}$ over ${\cH}_{\mathfrak P}$.
\item We have  $C(\alpha_k; I)={\cH}$ for each  $\alpha_k$
in a $\delta$-stair of $I$. So, we have  $U=Spec({\cH})$ and then
${\cH}[\partial]/I$ is ${\cH}$-flat.
\end{enumerate}
\end{Proof}

\begin{example} Let us denote $\CC[X]=\CC[x_1,\ldots,x_n]$ and consider the ideal of Example
\ref{ejemplo-bases-y-delta-bases}, i.e. $I=\left\{P_1,P_2\right\}
\subset \cD=\CC[X][\partial_1,\partial_2]$ where $$P_1=
x_1\partial_1+x_1\partial_2+b, \quad P_2=(x_2-x_1)\partial_2-d
 $$ with $b\in {\bf C}[x_1,x_3,\ldots,x_n]$ and $d\in
 \CC[x_3,\ldots,x_n]$. We will suppose $b$ is a multiple of $x_1$.
 In particular, $\cD/I$ is not a flat
 $\CC[x_1,\ldots,x_n]$--module, because the class of
 $\partial_1+\partial_2+b/x_1$ mod. $I$ has $x_1$-torsion.

We know by Example \ref{ejemplo-bases-y-delta-bases} and Remark
\ref{remark-after-sun} that $\{P_1,P_2\}$ is a Gr\"obner
$\delta$-base of $I$. Then a $\delta$-stair of $I$ is
$\left\{\exp^{\delta}(P_1),\exp^{\delta}(P_2)\right\}$, i.e.
$\left\{(1,0),(0,1)\right\}$.

Moreover, by Theorem \ref{tcdbR}, we have  $$ C((1,0);I)=C((1,0);
P_1,P_2)=\langle x_1\rangle, \qquad C((0,1);I)=C((0,1);
P_1,P_2)=\langle x_2-x_1\rangle.
 $$

Let us consider $J=C((1,0);I) C((0,1);I)$, i.e. $J=\langle
x_1(x_2-x_1)\rangle.$ By Theorem \ref{flateness},  ${{\bf
C}[X]}_{{\mathfrak P}}[\partial_1,\partial_2]/I_{{\mathfrak P}}$ is a
flat ${{\bf C}[X]}_{{\mathfrak P}}$-module for ${\mathfrak P}\in
U=Spec(\CC[X])\setminus V(J)$.
\end{example}

\begin{theorem}\label{finiteness}
Let $I$ be an ideal of ${\cD}$. The following are equivalent:
 \begin{enumerate}
 \item ${\cD}/I$  is a finitely generated ${\cH}$-module.
  \item For each $i=1,\cdots, n$ there exists $a_i\in \NN$ such that
   $$
   \alpha(i)=a_i\epsilon_i \in \Exp^\delta(I)$$ and
   $$C\left(\alpha(i);I\right)={\cH},$$ here $\epsilon_i$ is the
   $i$-th element of the canonical base of $\NN^n$.
 \end{enumerate}
 \end{theorem}

\begin{Proof}{.}
\noindent$1\Longrightarrow 2.$: For each  $i\in \{1,\cdots,n\}$ we
consider the  sub-${\cH}$-module  $M\subset {\cD}/I$ generated by the
set $$ \left\{1+I,
\partial_{i}+I,\cdots,
\partial^k_{i}+I, \cdots\right\}.$$ So, by the finiteness of $M$ over $\cH$, there exists
$\alpha({i})=a_i\epsilon_i\in \Exp^{\delta}(I)$ such that
$C\left(\alpha(i);I \right)={\cH}$, for some $a_i\in \NN$.

\noindent $2\Longrightarrow 1.$: Let us write $\overline \Delta=\NN^n
\setminus \bigcup_{i=1}^n \left(\alpha(i)+\NN^n\right)$. Let us
consider $M$ as the ${\cH}$-module generated by the finite set
$\left\{\partial^{\alpha}; \alpha \in \overline \Delta \right\}$. We
have ${\cD}=I+M$ and then ${\cD}/I$ is a quotient of $M$. Thus
${\cD}/I$ is finitely generated as ${\cH}$-module.
\end{Proof}


\begin{thebibliography}{90}
\bibitem[ASS-1]{Assi}{\sc A. Assi.} {\it Constructions effectives en alg\`{e}bre commutative.}
Tesis de Doctorado. Grenoble, 1991.

\bibitem[ASS-2]{assi1}{\sc A. Assi. } {\it On Flatness of generic projections.}
J. Symbolic Computation 18 (1994), 447-462.

\bibitem[ACG-1]{ACG1}{\sc A. Assi, F. Castro-Jim\'enez and J.-M. Granger. }
{\it The Gr\"obner fan of a $A_n$-module. } Journal of Pure and
Applied Algebra, 150 (2000), 27-39.

 \bibitem[ACG-2]{ACG2}{\sc A. Assi, F. Castro-Jim\'enez and J.-M. Granger. }
 {\it The standard fan of an analytic ${\cD}$-module. }
 Prepublicaci\'on n$^{\mbox{\underline o}} 47$
 Universidad de Sevilla, 1999.



\bibitem[BUCH]{buchberger1}{\sc B. Buchberger. }
{\it A theoretical basis for the reduction of polynomial to canonical
forms. } ACM. SIGSAM. Bull. 39 (1976), 19-29.



\bibitem[CAS-1]{castro-thesis}{\sc F.J. Castro. }{\it Th\'{e}or\`{e}me de division pour les op\'{e}rateurs
differentiels et calcul des multiplicit\'{e}s}. PhD thesis, Univ. Paris
VII, (Oct-1984).


\bibitem[CAS-2]{castro2}{\sc F.J. Castro. }{\it Calculs effectifs pour les id\'eaux d'op\'erateurs
diff\'erentiels. } In Travaux en Cours. G\'eom\'etrie Alg\'ebrique et
Applications, Tome III ,  Hermann, Paris, 1987. pp. 1-19.

\bibitem[CLO]{cox}{\sc D. Cox, J. Little and D. O'Shea. }{\it Ideals, varieties and
algorithms. An introduction to computational algebraic geometry and
commutative algebra. }Springer-Verlag, 1992.


\bibitem[IN-PA]{insa-pauer}{\sc M. Insa and F. Pauer.}
{\it Gr\"{o}bner bases in rings of differential operators. }In Gr\"obner
Bases and Applications. London Math. Soc. L.N.S. 251. Cambridge
University Press, Cambridge, 1998. 367-380.

\bibitem[MOR]{moreno-tesis}{\sc M.A. Moreno Fr\'{\i}as. } {\it M\'etodos
computacionales en los sistemas de ecuaciones en derivadas parciales.
} PhD. thesis. Universidad de Sevilla. Mayo 2000.


\bibitem[SAB]{sabbah-1} {\sc C. Sabbah.} {\it
Proximit\'{e} \'{e}vanescente. I. La structure polaire d'un ${\mathcal
D}$-module. (Appendice an collaboration avec F. Castro)}. Compositio
Math. 62 (1987), no. 3, 283--328.

\bibitem[TRI]{trinks} {\sc W. Trinks.}
{\it \"Uber Buchbergers Verfahren, Systeme algebraicher Gleichungen
zu l\"osen.} J. Number Theory, 10 (1978), no. 4, 475--488.

\end{thebibliography}
\end{document}